# Scale invariance versus translation variance in Nash bargaining problem

## Abstract

Nash's solution in his celebrated article on the bargaining problem calling for maximization of product of marginal utilities is revisited; a different line of argument supporting such a solution is suggested by straightforward or more direct reasoning, and a conjecture is raised which purports uniqueness of algorithm, namely his solution. Other alternative inferior algorithms are also suggested. It is argued in this article that the scale invariance principle for utility functions should and could be applied here, namely that utility rescaling u'=a*u is allowed, while translations, adding a constant to utility functions u'=u+b could not be applied here, since it is not invariant and leads to contradictory behavior. Finally, special situations of ownership and utilities, where trading is predicted not to take place at all because none is profitable are examined, and then shown to be consistent with the scale invariance principle.

---------------------------------------------------------------------------------------------------------------

## [1] Two-person Bargaining Problem

John Nash's article "The Bargaining Problem" relates to the classical economic problem of two parties freely bargaining with each other. In this situation, subsets of two original sets of items belonging to two individuals are possibly exchanged voluntarily if the sum of utilities is increased for both players.

The two opponents or players are each in possession of items that have utilities to both, and who engage in direct bartering without the use of money. Each wishes to convince the other to give away as many and most valuable items as possible in exchange for the fewest and least valuable items. A **proposal** of exchange involves a specific offer as well as a specific demand, where '**offer**' and '**demand**' refer to the details of the proposal. For example: "I offer the set R and I demand in exchange the set K".

We shall use the following notation for the utility functions of X and Y: *the-person-possessing-the item***ITEM***the-person-for-whom-utility-is-considered*. For example: $_x\mathbf{A}_y$ would mean the utility to Y of an item named A which X possesses originally.

Each possible exchange, whether reasonable or not, can be drawn on the $U_x$ vs. $U_y$ plane, where $U_x$ and $U_y$ refer to the marginal utilities of X and Y respectively, that is, the sum of utilities after an exchange is made minus the original sum of utilities of the items they owned, or equivalently, the sum of utilities of new items gotten (gained) minus the sum of utilities of items given out (lost).

It is very clear that only points in the first quadrant on the boundary/periphery (and not on the axes) are candidates for a solution, because any points inside (where there is room to go upward and/or to the right) are worse off for both players as compared with the boundary/periphery. For example, (3, 6) would not be considered at all if (5, 7) is also available with a different exchange, since the latter point offers more for X as well as more for Y. The existence of (3, 7) would also preclude the inferior point (3, 6) even though it's only worse off for Y. Moreover, Nash introduces the possibility of X and Y tossing a variable bias coin, with probability p on (0, 1) to decide between any two points on the periphery, hence creating a continuous line there. Note that this line includes linear connections between <u>all</u> points, and not only between adjacent ones on the periphery, thus resulting in an overall curve that is <u>not</u> the same as simply connecting adjacent points on the periphery, but rather a curve that encloses more area. The only question that remains then is "where exactly on that boundary?" and Nash's answer is that we look for an algorithm that maximizes $U_x * U_y$. Without any loss of generality, in this article, we would deemphasize or downplay the continuous line approach taken by Nash and instead focus only on those discrete points of actual deterministic exchanges. One severe consequence or limitation though of the discrete approach is that there are cases where maximum $U_x * U_y$ is not unique, whereas the continuous approach always guarantees uniqueness; hence an essential part of Nash's edifice was the joining of the points with this crooked yet continuous line. Examples of multiple maximum $U_x * U_y$ with the discrete approach abound.

# [2] Scale Invariance Principle

It is noted that Nash's solution point is scale invariant. The algorithm of maximizing $U_x * U_y$ would always point to the same transformed point (and the same items being exchanged – the same deal) regardless of scale. This is so because the maximum of a set of real numbers is scale-invariant in the sense that the same point (original/transformed) always serves as the maximum.

**IF**     $Max\{X_1, X_2, X_3, \ldots, X_n\} = X_i$  (for some unique i in the index to n )
**THEN**   $Max\{Q*X_1, Q*X_2, Q*X_3, \ldots, Q*X_n\} = Q*X_i$  (same i as above!!)
(where Q is **any** real positive number).

For example:

**IF**     $Max\{4, 7, 20, 3, 6, 10\} =$ **20**  (and it's the 3$^{rd}$ element)
**THEN**   Max **2**$*\{4, 7, 20, 3, 6, 10\} = Max\{8, 14, 40, 6, 12, 20\} =$ **40**
                                    (and it's also the 3$^{rd}$ element)

Moreover, **any** exchange where both X and Y are happy under original utility functions would continue to be rendered profitable to both after **any** transformation of scales **!**

To illustrate this, assume Y originally agrees to give X {A, B, C, D, … (M items)} and X agrees to give Y {P, Q, R, S, … (N items)} in exchange. To express the original agreement and satisfaction for both sides regarding this exchange, and taking first the point of view of X, we then must have:

   X extra satisfaction in obtaining {A, B, C, D, … M items}  >
                  >  X losses in not possessing {P, Q, R, S, … (N items)}
                               or :

     yAx+yBx+yCx +… M items  >   xPx+xQx+… N items

Now under scale transformation of the utility function of X by b (when b > 0), that is u'---> b*u we have to compare the two quantities

[ b*(yAx)+b*(yBx)+b*(yCx)+ … M items ] **and** [ b*(xPx)+b*(xQx)+… N items ]

   or, by pulling out the common factor b :

[ b*(yAx+yBx+yCx+ … M items) ]   **and**   [ b*(xPx+xQx+ … N items)]

for which the previous inequality above guarantees that the left side is bigger than the right side, regardless of the value of b. Same argument would apply to Y of course, and which completes the proof. The exact same line of argument above can be used to show that for any exchange refused by either X or Y or both under original utility function would also be resisted under any rescaling where b>0. To do that we simply flip the direction of the inequality above from > to < and all else follow just the same. **Hence any decision to trade for this two-person bargaining problem case is unaffected by any scale transformation, and therefore scale-invariant.**

The scale-invariance property of Nash's solution is not incidental. Solution (as well as everything else here) should not change if we re-scale utilities of X or Y or both. This is so because utility function itself is defined only up to scale. Hence if u is a utility function of X, than a*u is also a valid utility function. This expresses the fact that utility function is nothing but an expression of comparisons of preferences in owning one thing as oppose to another. For example, u(bike)=8 and u(book)=2 mean that the person derives four times as much satisfaction from a bike as compared with a book. Yet there is nothing intrinsic or important about the values 8 and 2. If we rescale everything by say 5, now u(bike)=40 and u(book)=10, we would still have the same ratio of satisfaction, namely four-fold, and it's this ratio that constitutes the real measure of preferences of utility, expressing the consistent fact that for that person four books satisfy him or her just as much as only one bike would, regardless of scale. Hence scale-invariance comes very naturally here. It is in this sense perhaps that utility values are thought of as using some arbitrary "satisfaction scale" or "happiness units" , just as the kilos, grams, and pounds are used for weights, especially or conveniently so when the least-valued item could be given the arbitrary value of 1.

## [3] Translations are not invariant!

Difficulties arise though when translations are considered. The other axiom of utility functions is that translated utilities are all equivalent expressions of satisfaction. In other words, that for any utility function u, u+b is also a valid utility function. But translations are harder to interpret. Say u(bike)=8 and u(book)=2, and now u is translated by +5, hence u(bike)=13 and u(book)=7, and so original ratios of satisfaction are totally upset! Moreover, if translation is performed using a value much larger than the largest value of the original utility function then in the limit each item is assigned equal value – clearly a false representation of preferences and satisfaction order! For example, in the case above where u(bike)=8 and u(book)=2, and where b is 1000 say, translated utilities are u(bike)=1008 and u(book)=1002, and are of almost identical values.

**Moreover, here trades are NOT even translation-invariant! In other words, a given exchange acceptable to both X and Y under original utility function would be rendered as totally unacceptable under a certain utility translation!**

To illustrate this, assume Y originally agrees to give X {A, B, C, D, … (M items)} and X agrees to give Y {P, Q, R, S, … (N items)} in exchange. To express the original agreement and satisfaction for both sides regarding this exchange, and taking first the point of view of X, we then must have:

X extra satisfaction in obtaining {A, B, C, D, …  (M items)}   >
> X losses in not possessing {P, Q, R, S,  … (N items)}

or :

yAx+yBx+yCx +… M items  >   xPx+xQx+… N items

Now under translation of utility function of X by b, that is u ----> u+b we have to compare the two quantities

[(yAx+b)+(yBx+b)+(yCx+b)+… M items ] **and** [ (xPx+b)+(xQx+b)+… N items ]

or :

[(yAx+yBx+yCx +… M items) +M*b ]   **and**   [ (xPx+xQx + … N items) + N*b ]

for which the previous inequality above does **<u>not</u>** by any means guarantee a similar inequality here, unless b is relatively very small. If M=N or if M>N then the inequality still holds true, but whenever M<N things depend on the relative value of b. When b is extremely large, satisfaction in general (in the limit) is calculated by simply counting the number of items exchanged, regardless of type, and in such a case all depends solely on the relative value of N and M!

Another strong consideration against allowing translation of utility functions is for situations when an item deemed totally useless is being possessed, say children clothes for a family without kids, etc., and in which case the original utility value for such an item is naturally zero. Yet, in this case, translation would result in a non-zero positive value for such an item, while rescaling would always guarantee a consistent zero value no matter what scale is chosen.

**Hence it is necessary to disallow translation of utility function altogether and insist on rescaling as the only possibility at least in the context of this particular two-person bargaining situation.**

It is for this property of the utility function leading to scale-invariance here that neither X nor Y would ever "complain about fairness" if solution supposedly is at $(X_0, Y_0)$ and $X_0 > Y_0$ or visa versa. This is so because with a simple rescaling we could easily reverse the situation and obtain $X_0 < Y_0$ and then perhaps X would complain now, etc. On the contrary! The correct perspective or view here runs much deeper! Neither would complain about fairness at all, as neither of them is concerned about how the other person is fairing, nor comparing his or her overall utility to that of the other. Rather, all Y is trying to do is to obtain a point as high as possible on the $U_x$ and $U_y$ plane, while X is striving to obtain a point as much to the right as possible. The only solution point that could consistently be considered unfair to say Y is $(X_0, 0)$, namely on the X axis, in which case no rescaling could ever reverse the favorable situation that X enjoys here.

## [4] Naturalness of Maximum Marginal Utilities

Deemphasize the continuous line approach taken by Nash and focusing solely on those discrete points, it is rather intuitive or natural that we'll look for or find the solution at maximum $U_x * U_y$ rather than anything else. We are in need of an algorithm that is **(1) fair to both, (2) beneficial to both, (3) unique, and (4) scale invariant.**

First, it is clear that **both** $U_x$ and $U_y$ must be somehow involved or present in any expression or algorithm, otherwise one would be considered as being discriminated against the other. Of the four basic arithmetic operations, only addition and multiplication can be used interchangeably without any preference regarding operands, while subtraction $U_x - U_y$ and division $U_x / U_y$ are intrinsically preferential operations, treating $U_x$ and $U_y$ in different ways, hence should be excluded here all together. Expressions such as $(U_x^5 * U_y^2)$ or $(U_x)^{U_y}$ would also not be acceptable as it discriminates against Y and favors X or visa versa. Secondly, what is needed here is a decisive algorithm, operation, or procedure that points to a unique point, not multiple ones. Let us now then compare maximizing $U_x+U_y$ (addition) with maximizing $U_x * U_y$ (multiplication). Evidently multiplication is superior here because product is one operation that is maximized only when both $U_x$ and $U_y$ are not small (and especially not close to zero), and so in that sense, by maximizing the product of the two marginal utilities we consider the interests of both X and Y equally and fairly. On the other hand, maximization of $U_x + U_y$ can be done (if we wish) totally at the expense of one of them. We might even find cases where the sum is being maximized while one quantity say $U_x$ is zero**,** a situation or a solution which does not reflect our case given that trading is voluntary for both players and that they are of equal knowledge of utilities and of equal bargaining power. Another strong consideration against maximizing $U_x + U_y$ is that it's not scale invariance. For example, if we rescale utility of (only) X by a huge factor, things would swing much in X favor as maximization of $U_x + U_y$ will now end up focusing much more on X than on Y. Another drawback of maximum $U_x + U_y$ is that it could yield multiple solutions whenever boundary points farthest from the origin (i.e. the best points) fall on the equidistance line $U_x + U_y = C$. And what other methods would X and Y would agree upon? Maximizing $(1/U_x)*(1/U_y)$? Minimizing $U_x * U_y$? Certainly not, as it would hurt both. Finding the point closest to the meaningless line X=Y which vary with scale? surely not! Both X and Y would like to maximize quantities here, not to minimize, nor to set quantities equal to anything when scales are totally meaningless. **There seems to be no other reasonable algorithm we can come up with besides maximizing the product $U_x * U_y$!**

Consider the case of 2-object gravitational force, $F = G*M_1*M_2/R^2$. Say we have a fixed amount of matter, 10 kilograms, where mass is flexible so that we can draw from one box of matter and into another as much as we wish, in the same vein as zero-sum or constant-sum situations. Distance is fixed at 10 meters say. The question then arises, under what arrangement of matter do we get maximum gravitational force? Here is the table:

| $M_1$ | $M_2$ | $M_1*M_2$ | F |
|---|---|---|---|
| 0 | 10 | 0 | 0 |
| 1 | 9 | 9 | 6.00 E-12 |
| 2 | 8 | 16 | 1.06 E-11 |
| 3 | 7 | 21 | 1.40 E-11 |
| 4 | 6 | 24 | 1.60 E-11 |
| **5** | **5** | **25** | **1.66 E-11** |
| 6 | 4 | 24 | 1.60 E-11 |
| 7 | 3 | 21 | 1.40 E-11 |
| 8 | 2 | 16 | 1.06 E-11 |
| 9 | 1 | 9 | 6.00 E-12 |
| 10 | 0 | 0 | 0 |

Hence, by keeping the two weights as even as possible, that is by dividing them equally, we obtain maximum product, while extreme inequality yielded zero. Hence, in general it can be said that maximizing product of constant-sum variables entails the most equal division possible! Yet this conclusion can not be separated from the fact that rescaling here say only $M_1$ leaving $M_2$ in tact is forbidden. Only simultaneous and equal rescaling is allowed here to uphold any meaning to the restriction $M_1 + M_2 = 10$, and that rescaling should be applied to the number 10 as well, comparing apples to apples not to bananas.

In our context, this would point to a "fair" solution at (C/2, C/2) for all the discrete points on the line of $U_x + U_y = C$. Surely, almost all typical curves in bargaining games are convex as opposed to straight lines, and one-sided rescaling is surely allowed, yet some principle of "fairness" in maximizing products is demonstrated by this physical example.

An intriguing and totally different line of thought is that mathematically speaking **there may not exist any other algorithm arriving at a unique point except maximization of $U_x*U_y$ - if we wish to be fair to both, to benefit both as much as possible, and for the algorithm to be scale invariant.**

In other words, there is simply no other algorithm that could differentiate between the various points available on the periphery and point to a unique one except maximizing $U_x*U_y$ ! If so, then say Y is unhappy about the end result of maximizing $U_x*U_y$, and so Y tries hard to convince X to choose another point $(X_0, Y_0)$. In that case Y wouldn't be able to back up this point by appealing to <u>any</u> algorithm. Y would simply have to arbitrarily insist or ask X to move to that point without being able to say anything in support of his suggestion!

## [5] Other Algorithms

The only other scale-invariant algorithm perhaps that could be proposed here, uniquely pointing to a special point on the periphery in a way which is (at first glance) fair to both X and Y is the following:

**Algorithm A)** Count all those discrete points on the periphery (excluding points on the axes) and choose the middle one if they are odd, or the midpoint (tossing a biased coin) between the two most-centralized points if they are even, in the same spirit as in the definition of the median. This algorithm is scale invariant!

This algorithm could work well for our discrete-points approach in this article. But with the inclusion of lines generated by tossing a deciding coin, - the continuous approach suggested by Nash - there would be nothing special about this (supposed) most-centralized point! Nor could we choose some midpoint, by way of measuring the distance, along the zigzagged line for the continuous approach, as it is scale-dependent, and hence meaningless.

Furthermore, the apparent fairness in this approach is misleading and the disadvantage is revealed in cases where the points on the periphery crowd out unevenly close to one axis and are sparse and spread out next to the other axis, for example:

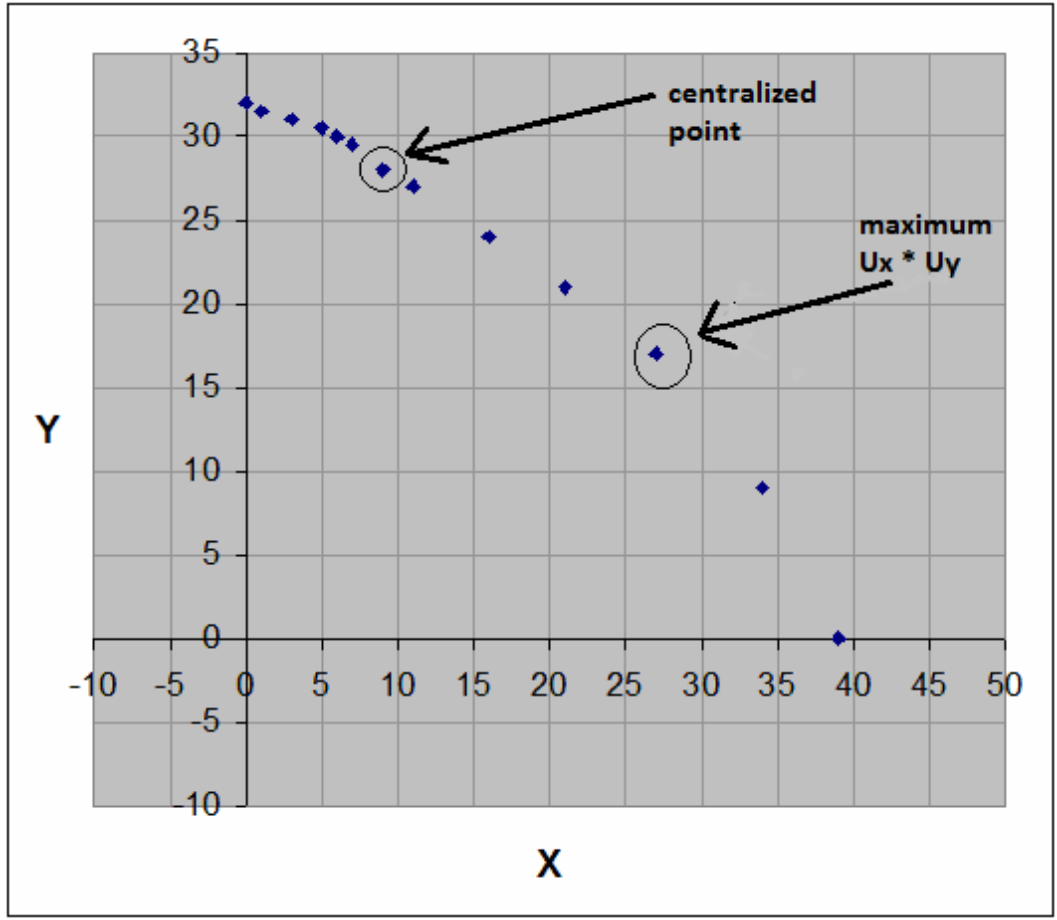

| Ux | Uy | product Ux*Uy | comments |
|---|---|---|---|
| 0 | 32.0 | 0.0 | excluded |
| 1 | 31.5 | 31.5 | |
| 3 | 31.0 | 93.0 | |
| 5 | 30.5 | 152.5 | |
| 6 | 30.0 | 180.0 | |
| 7 | 29.5 | 206.5 | |
| **9** | **28.0** | **252.0** | **centralized** |
| 11 | 27.0 | 297.0 | |
| 16 | 24.0 | 384.0 | |
| 21 | 21.0 | 441.0 | |
| 27 | 17.0 | **459.0** | **max Ux*Uy** |
| 34 | 9.0 | 306.0 | |
| 39 | 0.0 | 0.0 | excluded |

In this case, the central point (9, 28) is not quite acceptable to X who would very much prefer the point (27, 17) where $U_x*U_y$ is maximized. Note that this dislike of X here regarding choosing the most centralized point is independent of scale. Should X rescale his utility function by a factor much bigger than 1 or by a tiny fraction of 1 (close to 0) his or her dislike of such a solution would remain in tact and X would always seek to convince Y to agree on maximum $U_x*U_y$ - a proposition Y would not view negatively nor as unreasonable.

A very different approach combining the concern for the two issues at hand here - fairness and scale invariance, is to force one final 'fair' rescaling before performing any algorithm or calculations. This done, any subsequent algorithm can now be made fairly disregarding scale altogether, and permitting procedures that are scale-dependent. This is done by equating (via rescaling) the X and Y distances of the two extreme points regarding all 4 quadrants. An alternative way would be to equate X and Y distances of the curve considering only the first quadrant for the purpose, hence rescaling both X and Y such that (r, 0) and (0, s) [or their closest points if none of the points falls on the axes] would be of the same length after transformation, that is (k, 0) and (0, k). We shall call this "**equitable rescaling**".

**Algorithm B)**  Perform equitable rescaling, then choose the point which maximizes X+Y. The difficulty here is that there may be no unique solution here, even though the zigzagged line connecting the points is concave.

**Algorithm C)**  Perform equitable rescaling, then choose the point closest to the 45 degrees line X=Y.

**Algorithm D)**  Perform equitable rescaling, then measure the total distance on the zigzagged line to choose the point half way through on that curve. We could either do it by creating the zigzagged line considering strictly only adjacent points, or the more encompassing approach given by Nash.

## [6] Cases Having no Profitable Trades

The force driving the exchange is the different valuations placed by X and Y for some of the items, implying that each possible exchange could in principle represent a totally different pair of sums of utilities. There would have been no incentive whatsoever to trade had utilities were the same for both players for each item, that is if xZx = xZy and yZy = yZx for any item Z. In such a case, all the points on the Ux and Uy plane appear along the line Y = -X going through the origin, since in any exchange here the gain for X (net of some utilities gained and some lost) is exactly the loss for Y. For example, say X gives out three items with total utilities 30 (valued equally), and gets from Y two items with total utilities 34 (valued equally), resulting in Ux = +4 and Uy = -4.

Should we rescale here, say utilities of X are to be multiplied by a factor of 2; this would only mean that the X axis was stretched to double its original size, resulting in a line through the origin with a slope of -0.5 instead. It was proven earlier in this article, in absolute generality, that all trade decisions are independent of scale.

Another situation where no trade could ever be made possible is depicted as follow: for each item X possesses the satisfaction derived from it for X is greater than that for Y, AND ALSO visa versa. In other words, for any item Z possessed by either X or Y we have iZi>iZj (indexes i and j representing X or Y). This precludes any possibility of trading. For example:

| X possesses | radio | Laptop | Book | watch | Pen |
|---|---|---|---|---|---|
| utility to X: | 11 | 8 | 5 | 6 | 4 |
| utility to Y: | 4 | 3 | 1 | 1 | 3 |

| Y possesses | bike | TV | Cell | chair |
|---|---|---|---|---|
| utility to Y: | 7 | 3 | 11 | 6 |
| utility to X: | 4 | 2 | 10 | 5 |

Let's see why any possibility of one-to-one trade is precluded. Repeated use of the transitive rule for inequality will be made here, namely that if A>B and B>C then A>C. Say X possesses item A, and Y possesses item C, among some other items, and our particular situation is depicted as xAx>xAy and yCy>yCx, hence a trade of item A for item C in which X is benefiting from the exchange gives yCx>xAx (gaining from C more than losing on A) and this implies that yCx>xAy (because of xAx>xAy). Now, from yCy>yCx we see that yCy>xAy, meaning that the loss of C for Y is greater than Y's gain in acquiring A, and therefore Y could not be benefiting from such a trade!

Similar argument can be made in general for ANY numbers of items traded, say Y gives X {A, B, C, D, ... (M items)} and X gives Y {P, Q, R, S, ... (N items)}, keeping in mind that in here the condition implies that for any item Z we have iZi>iZj - namely that the original owner obtains more satisfaction than the other person for any item possessed.

Assume X would benefit and agree to such a trade as above, this would imply that (gains>losses):

[yAx+yBx+yCx +... M items]  >  [xPx+xQx+... N items]

and employing rule iZi>iZj to the right side of the inequality we get:

[yAx+yBx+yCx +... M items]  >  [xPy+xQy+... N items]

Now comparing the left side of this inequality to [yAy+yBy+yCy +... M items], we note that the latter quantity is superior (again because of iZi>iZj), and so we finally get:

[yAy+yBy+yCy +... M items] > [xPy+xQy+... N items]

Hence clearly the loss of Y from not owning {A, B, C, etc.} is greater than Y's gain from newly acquired items {P, Q, R, etc.}, and therefore Y would refuse to trade as such!

The scale invariance principle applies. Doubling utilities of X would make the gap with Y wider for items X possesses but narrow it or even reverse the situation for items Y possesses. If X or Y refuses any trade under some given utility functions, then they should continue to resist such an exchange under any rescaling of utilities as was proven earlier.

A third example of a situation where no trade could ever be made possible is when the SUM of utility values for (say) X regarding ALL the items possessed by Y is LESS than the value X places on ANY of his own items (or equivalently: LESS than the MIN value X has for his own items.) In this case no amount of satisfaction from items belonging to Y can ever be awarded to X to compensate for the loss of even one single item he possesses, even the one least valued.  For example:

| X possesses: | Radio | Laptop | book | watch | Pen | | |
|---|---|---|---|---|---|---|---|
| utility to X: | 14 | 12 | 13 | 12 | 13 | | |
| utility to Y: | 10 | 16 | 7 | 17 | 5 | | |

| Y possesses: | Bike | TV | cell | chair | | | |
|---|---|---|---|---|---|---|---|
| utility to Y: | 2 | 8 | 13 | 2 | | | |
| utility to X: | 3 | 1 | 4 | 3 | | = | 11 |

Any rescaling of utility values of X would involve the same multiplicative factor G say, for the transformation for items possessed by X as well as for items possessed by Y, hence transformed utility sum of all the items Y possesses is still for X less than the transformed value for any of his own items. This is consistent with the general scale invariance principle proven earlier.

It is noted that in this last example, a one-way inequality is a sufficient condition to preclude trading, as oppose to the earlier example where mutual inequality was a necessary condition, namely iZi>iZj  as well as  jZj>jZi .

## [7] General Considerations

To facilitate numerical analysis Nash makes the assumption that 'ultimate' or 'ideal' utility for each player, his or her overall goal, is simply the <u>sum</u> of the utilities of the individual items he or she possesses.  In reality this is certainly not true.  The first loaf of bread has much higher overall utility than the 10th loaf.  One house and a loaf of bread have a much higher overall utility than either two loafs (full belly but staying outside in the cold) or two houses (spacious living space but with a hungry stomach). Nonetheless, an analysis based on the more realistic assumption of the combinatorial effect on overall utility would follow the same line of reasoning as outlined in Nash's article given the more detailed 'ultimate utility function' which yields unique values for each and every combination of items, leading to a more realistic set of points on the $U_x$ and $U_y$ plane for determinations of max $U_x*U_y$ .

It is interesting to note that the uniqueness of solution referred to in Nash's article is not about the expected details of the exchange (what items exactly are being exchanged) but rather about the expected overall utilities for the two players after a deal is struck, values which at times could be achieved in more than one way of exchange and in which case exact details of exchange couldn't be predicted at all. To demonstrate this, imagine the case where player X possesses two items, A and B, and derives equal satisfaction from either one, and where Y also derives equal satisfaction here, that is, $xAx = xBx$ and $xAy = xBy$, while solution involves X giving out to Y one of these two items. There is no way to predict which item exactly X will give up, so solution is not unique in that sense.

Implicitly, Nash offers the use of the brute force method, namely, to check the entire list of all possible exchanges and sort out maximum Ux*Uy. Even with the elimination of negative and zero marginal utilities this leaves the computer performing combinatorial number of operations which increases exponentially as the number of items possessed becomes large. One intriguing question is whether or not there exists a way to cut directly to the solution, in the same vein as when the gradient seeks out the maximum and minimum of multivariable functions by cutting right through to it, avoiding the brute force tedious way. Attempts to devise an algorithm consisting of numerous little logical or very "reasonable" exchanges leading to the final solution are doom to failure. A sequence of one-for-one exchanges of items that increase sums of utilities to both players will not do. This is so because such an algorithm couldn't capture the benefits presented by simultaneous multi exchanges.

It is hard to imagine how two rational players during Adam Smith's era living without the benefit of computers could calculate and arrive at the optimal solution when number of all possible trades is really huge as in the case when items to be exchanged are quite numerous. If X possesses P items and Y possesses Q items, then there are $2^{(P+Q)}$ possible exchanges, including the case of the null exchange - where nobody gives anything out nor takes anything, that is, keeping the original situation of ownership in tact. Many of those $2^{(P+Q)}$ possible exchanges though converge onto identical points on the $U_x$ and $U_y$ plane whenever resultant $U_x$ and $U_y$ are equivalent, hence total number of points typically is much less than that and often folds into something like a two-thirds, a half, or even around a third of the maximum potential of $2^{(P+Q)}$ points.

# References

John Nash J.F. (1950), The Bargaining Problem, Econometrica 18, 155-162.

Alex Ely Kossovsky
akossovs@yahoo.com